\documentclass{article}
\usepackage{poly}

\hypersetup{
colorlinks=true,
citecolor=black,
linkcolor=black,
anchorcolor=black
,filecolor=black,
menucolor=black,
urlcolor=black,
pdftitle={Polyhedral approximations of Riemannian manifolds},
pdfauthor={Anton Petrunin}
}

\begin{document}

\title{Polyhedral approximations\\of Riemannian manifolds}
\author{Anton Petrunin}
\date{}
\maketitle

\begin{abstract} 
We give a condition on the curvature tensors of Riemannian manifolds that admit Lipschitz approximation by polyhedral metrics with curvature bounded below or above.

We show that this condition is also sufficient for the existence of local approximations.
We conjecture that it is also sufficient for the global approximations and prove it in some special cases.
\end{abstract}

\section*{Introduction}

\parbf{Positive cosectional curvature.}
Let $\EE^m$ denotes  $m$-dimensional Euclidean space.
Denote by $\A^4(\EE^m)$ the space of all curvature tensors of Riemannian manifolds on the tangent space $\EE^m$.
The $\O(m)$-rotations of $\EE^m$ induce isometric rotations of space of tensors in $\A^4(\EE^m)$.
The subset of $\A^4(\EE^m)$ will be called \emph{$\O(m)$-invariant} if it is invariant with respect to these rotations.

Let $Q_{\SS^2\times \RR^{m-2}}$ be the curvature tensor of $\SS^2\times \RR^{m-2}$ with the standard metric.
Denote by $\mathcal{S}^*$ the minimal convex $\O(n)$-invariant cone in $\A^4(\EE^m)$ that contains $Q_{\SS^2\times \RR^{m-2}}$.

Let $M$ be a Riemannian manifold and $p\in M$;
denote by $\Rm_p$ the curvature tensor of $M$ at $p$.
We say that \emph{cosectional curvature} of $M$ at $p$ is at least (at most) $\kappa$
if \[\Rm_p-\kappa\cdot Q_{\SS^m}\in\mathcal{S}^* \quad\text{or, respectively,}\quad-\Rm_p-\kappa \cdot Q_{\SS^m}\in\mathcal{S}^*,\]
where $Q_{\SS^m}$ denotes the curvature tensor of unit $m$ sphere $\SS^m$.
Briefly, these conditions can be written as $\cosec(\Rm_p)\ge \kappa$ and $\cosec_p\ge \kappa$ or respectively $\cosec(\Rm_p)\le \kappa$ and $\cosec_p\le \kappa$.
The reason for the name \emph{cosectional} will become clear in 1.E.

We will write $\cosec M\ge \kappa$ (or $\cosec M\le \kappa$) if $\cosec_p\ge \kappa$ (respectively, $\cosec_p\le \kappa$) for any point $p\in M$.

The following theorem gives an if-and-only-if condition for the existence of local approximations of a Riemannian manifold by polyhedral spaces with a lower curvature bound.

Recall that a length-metric space is called \emph{$\kappa$-polyhedral} if it admits a triangulation such that each simplex is isometric to a simplex in the model space of curvature $\kappa$.
The $0$-polyhedral spaces are called \emph{Euclidean polyhedral spaces},
the $1$-polyhedral spaces are called \emph{spherical polyhedral spaces},
and
the $(-1)$-polyhedral spaces are called \emph{hyperbolic polyhedral spaces}.

\begin{thm}{Local Theorem}\label{0.1}
Let $P_n$ be a sequence of $m$-dimensional $\kappa$-polyhedral spaces
with curvature $\ge \kappa$ in the sense of Alexandrov.
Assume $P_n$ Lipschitz converges to a Riemannian manifold $M$, then $\cosec M\ge \kappa$.

Moreover, if $M$ is a Riemannian manifold with $\cosec M\ge\kappa+\eps$ for some $\eps>0$,
then each point of $M$ has a neighborhood that is a Lipschitz limit of a sequence of
polyhedral spaces with curvature $\ge \kappa$.
\end{thm}

\begin{thm}{Global theorem}\label{0.2}
Let $M$ be Riemannian $m$-manifold with $\cosec M\z\ge \kappa+\eps$ for some $\eps>0$.
Assume that $M$ (or its finite cover) has a stably trivial tangent bundle.
Then $M$ admits a Lipschitz approximation
by a sequence of $m$-dimensional $\kappa$-polyhedral metrics with curvature
$\ge \kappa$.
\end{thm}

We expect that the condition on the tangent bundle is not essential.
The simplest interesting example that is not covered by the theorem is $\CP^2$ ---
the complex projective plane with the canonical metric.
This problem is discussed further in Section \ref{sec:Problem section}.

The global theorem can be reduced to the cases $\kappa= \{-1,0,1\}$.
In the first and last cases using rescaling one can get an approximation
of $(M,g)$ with polyhedral metrics with curvature $\ge \mp1$.
In case $\kappa=1$, the condition $\cosec M\ge 1$ implies, in particular, that the curvature operator of $M$ is strictly positive (see 1.E).
In particular, the Micallef--Moore theorem \cite{MM} implies that universal cover $\widetilde M$ of $M$ has to be homeomorphic (and by B\"ohm--Wilking theorem \cite{BW} diffeomorphic) to
the standard sphere.
In particular, $\widetilde M$ has stably trivial tangent bundle.
Therefore we get the following:

\begin{thm}{Corollary}\label{0.3}
A complete Riemannian $m$-dimensional manifold $M$ admits a Lipschitz approximation by
spherical polyhedral metrics with curvature $\ge 1$ 
if and only if $\cosec M\ge 1$.
\end{thm}

\parbf{About the proofs.}
The necessity of the curvature bound follows from the fact that all curvature
of a polyhedral metric stays on hyperedges;
that is, the simplexes of codimension $2$ of the triangulation of the polyhedral space.
Around every hyperedge, the metric looks like $C\times \RR^{m-2}$,
where $C$ is a two-dimensional cone.
Thus the curvature at the vertex of $C$
looks like the curvature of $\SS^2$ with zero radius,
 and this allows us to view the curvature at the edge as the curvature of
$\SS^2\times \RR^{m-2}$ (which is $(x\wedge y)^2$) multiplied by a
Hausdorff measure of edge.
If $M$ is approximated by polyhedral metrics, then the curvature
tensors of different edges could mix with each other;
that is, the limit manifold $M$ must have curvature tensor that is a convex combination of the curvatures of the above form;
in other words, $\cosec M\ge 0$.

It gives only an idea; the real proof is very technical.
This part of the proof is not included in this paper.
In smooth settings, this work was recently done by Nina Lebedeva and the author \cite{LP}.
The proofs in polyhedral settings are essentially the same. Both parts are based on Perelman's unpublished result that is discussed in Appendix \ref{appendix:B}.

To prove the sufficient condition, we construct an isometric embedding of $(M^m,g)$ into $\RR^q$ such that its image is an intersection of $q-m$ locally convex hypersurfaces with obtuse angles between their outer normals at any point of~$M$.
This condition on the angles alone implies $\cosec M>0$.

Next, we consider an approximation of convex hypersurfaces
by convex polyhedral hypersurfaces with the same condition on the angles between the outer normals.
The needed polyhedral approximation is the intersection of these
polyhedral hypersurfaces. 
That proves the local theorem.

Now let us describe the idea for the global theorem.
For simplicity assume that $M$ is simply-connected.
Using that $\T M$ is stably trivial we realize $M$ as an intersection of open convex hypersurfaces with the same conditions on angles, and then do the same approximation as above.
The proof of this last representation is technical. 
Note that if such representation exists, then $\N M$, the normal bundle of
$M$, is trivial; in particular, the tangent bundle $\T M$ should
be stably trivial.
The latter explains how the condition on the tangent bundle in the global theorem comes into the game.

\parbf{Acknowledgments.}
The main part of this work was done during my stay at the IHES in 1999--2000.
I would like to thank this institute for  its support and hospitality.
I want to thank Grigory Perelman for sharing ideas and making me interested in this problem long before this publication; 
Vladimir Voevodsky for bringing the paper
of David Hilbert to my attention. 
I want to express my very special thanks to Jost Eschenburg and Sergei Kozlov who constructed for me
necessary examples of curvature tensors and pulled me out of a
dead-end in this research.
I want to thank Rostislav Matveyev, Dmitri Panov, and Thomas Sharpe for their helpful and interesting conversations and letters.

\section{Notation, definitions, and preliminaries}

\parbf{1.A. Polyhedral spaces with curvature bounded below or above.}

A connected simplicial complex $P$ is called \emph{pseudomanifold}
if the link of each simplex in $P$ is connected or $\SS^0=\{-1,1\}$.

Let us denote by $\MM^q[\kappa]$ the $q$-dimensional simply connected space of constant curvature $\kappa$.
A pseudomanifold $P$ equipped with a metric 
such that each
simplex is isometric
to a simplex in $\MM^m[\kappa]$
is called
\emph{$\kappa$-polyhedral space}. 

A $\kappa$-polyhedral space has curvature bounded below  if  the sum of angles
around any hyperedge (that is, simplex of codimension $=2$) is $\le 2\pi$.
In this case, the polyhedral space has curvature  $\ge \kappa$ in the sense of Alexandrov; see \cite{AKP}.

In all that follows we will assume $\kappa=0$,
but if it is not specially
mentioned everything below is true for any $\kappa$;
one only has to exchange $\RR^q$ to $\MM^q[\kappa]$.


\parbf{1.B. Convex submanifolds of higher codimension.}

\begin{thm}{Definition}
 A submanifold $M \i \RR^q$ is called
\emph{locally convex} if for each point $p\in M$ there is a collection of strictly convex hypersurfaces $F_i$,
such that $U=\bigcap_i F_i$ is a neighborhood of $p$ in $M$ and at each point of $U$ the angle between
outward normals to any pair of $F_i$ is obtuse.
\end{thm}

Note that we do not assume that the hypersurfaces $F_i$ are closed subsets in~$\RR^q$.

If $M$ is $C^2$-smooth, then the above property is equivalent to the following condition:
 at each point $x\in M$ there is an orthonormal
basis $\{e_i\}\i \N_x M$, where $\N M$
is the normal bundle of $M\i \RR^q$, such that the representation of
the second fundamental form
$s\:\S^2(\T_xM)\to \N_x M$ in this basis $s=\sum_i e_i\cdot s_i$
has all positive-definite forms $s_i\in \S^2(\T)$.

The latter property can also be reformulated the following way: a smooth submanifold is
locally convex if it can be viewed locally as a convex hypersurface
in a convex hypersurface in $\dots$ in $\RR^q$.

It is easy to see that locally convex submanifolds have positive curvature (for the induced intrinsic metric in the sense of Alexandrov).
If the submanifold is smooth, then one can say  more about its curvature tensor;
the latter is done in the following subsection.

\parbf{1.C. Curvature tensors of submanifolds.}
Here we introduce an extrinsic curvature for submanifolds.
This subsection is based on \cite[3.1.5]{Grom-PDR}.

Let $\T$ be a vector space with a scalar product, $\T^n$ will denote
its tensor power of degree $n$, $\S^n(\T)$ and $\Lambda^n(\T)$ will
denote respectively subspace of symmetric and antisymmetric
elements of $\T^n$. The scalar product on $\T$ induces a
scalar product on $\T^n$ and all its subspaces.

The following subspace of $\T^4$
$$\A^4(\T)=\Lambda^4(\T)^\bot\cap \S^2(\Lambda^2(\T))$$
is formed by all possible curvature tensors
on the tangent space $\T$.
In an equivalent way this subspace $\A^4\i \S^2(\Lambda^2(\T))$
can be described  as the space of all tensors in $\S^2(\Lambda^2(\T))$
satisfying the first Bianchi identity
$$\Rm(X,Y,Z,W)+\Rm(Y,Z,X,W)+\Rm(Z,X,Y,W)=0.$$
In particular, the subspace $\A^4(\T)$ does not depend on the choice of scalar product on $\T$.

\medskip

Let $M\i \RR^q$ be a smooth submanifold and $s_x\:\S^2(\T_xM)\to \N_x M$ its second fundamental form at $x\in M$;
here $\T M$ and $\N  M$ are tangent and normal bundles over $M$ respectively.
Consider the \emph{$\Phi$-curvature tensor}
$$\Phi (X,Y,Z,W)= \langle s(X,Y),s(Z,W)\rangle,$$
here $\Phi$ is a section of $\S^2(\S^2(\T M))$.

Tensor $\Phi$ can be written as
$$\Phi(X,Y,Z,W)=E(X,Y,Z,W)+\tfrac 1 3\cdot(\Rm(X,Z,Y,W)+\Rm(X,W,Y,Z))
\eqno(*)$$ 
where $E$ is the total symmetrization of $\Phi$; that is,
$$E(X,Y,Z,W)=\tfrac 1 3\cdot
(\Phi(X,Y,Z,W)+\Phi(Y,Z,X,W)+\Phi(Z,X,Y,W))\in \S^4(\T),$$
and
$$\Rm(X,Y,Z,W)=\Phi(X,Z,Y,W)-\Phi(X,W,Y,Z)\in \A^4(\T)$$
is the Riemannian curvature tensor of $M$.

Tensor $E$ represents the \emph{extrinsic curvature} of $M$.
Note that $E\in \S^4(\T)\z\i \S^2(\S^2(\T))$.
Tensor $E$ measures wrinkling of the embedding --- the more it is wrinkled the bigger $E$ gets.
Note that
$f(X)=E(X,X,X,X)=|s(X,X)|^2$
is homogeneous polynomial of degree $4$ on $\T$ 
and it describes $E$ completely.
Therefore in some sense, the $E$-tensor is a higher order analog of the metric tensor.

There are two reasons to use tensors $\Phi$ and $E$: 
\begin{itemize}
\item The $\Phi$-curvature depends only on elements of $\T$;
in particular, it does not depend on the dimension of the ambient space.
This way we can study isometric embeddings without direct referring to the ambient space.
\item Direct
construction shows that $\Phi$ describes the second fundamental
form $s\:\S^2(\T)\to \N$
up to an isometric rotation of $\N$;
that is, two second
fundamental forms $s_1,s_2\:\S^2(\T)\to \N$ 
give the same tensor $\Phi\in\S^2(\S^2(\T))$ if and only if there is an isometric rotation
$j\:\N\to \N$, such that $j\circ s_1=s_2$. 
In particular, since
$\Phi$ is a sum of Riemannian curvature tensors and $E$, we have
that if $(M,g)$ is a Riemannian manifold and $(M,g)\to \RR^q$ is
an isometric embedding, then $E$-tensor together with $g$ describes
the second fundamental form at each point up to an isometric
rotation.
\end{itemize}

\parbf{1.D. Positiveness in $\bm{\S^2(\S^2(\T))}$
and convexity of submanifolds.}
Most of this subsection is extracted from \cite[2.4.9B(4)]{Grom-PDR}.

Given an open convex cone $\mathcal{C}$
in a Euclidean space $\EE^n$, set
$$\mathcal{C}^*=\set{r\in \EE^n}{\langle r,r'\rangle >0 \text{ for all }r'\in \mathcal{C}}.$$

\begin{thm}{Definition} A tensor $\Phi\in \S^2(\S^2(\T))$
is positive ($\Phi> 0$), if there is a
representation $\Phi=\sum_i s_i^2$, where $s_i$ are positive-definite forms on~$\T$.
If $i\:M\to \RR^q$ is a smooth embedding we will write
$\Phi(i)>0$ if the $\Phi$-tensor of $i(M)\i \RR^q$ is positive
at $i(x)$ for all $x\in M$.
\end{thm}

The cone of positive tensors in $\S^2(\S^2(\T))$ forms a convex
$\GL(\T)$-invariant cone of tensors;
that is, this cone is invariant with respect to the action of $\GL(\T)$ induced on $\S^2(\S^2(\T))$.
If $\dim \T\ge 2$, then
there are other $\GL(\T)$ invariant cones in $\S^2(\S^2(\T))$, one of
these cones will be of particular interest; it is the cone
of all elements $\Phi=\sum_i s_i^2$ for arbitrary elements $s_i\in
\S^2(\T)$. This cone describes all elements of $\S^2(\S^2(\T))$ that
can appear as $\Phi$-curvature of a submanifold.

Note that existence of representation of the second fundamental form
$s\z=\sum_i s_i\cdot e_i$ with positive-definite $s_i\in \S^2(\T)$
implies, in particular, that $\Phi\z=\sum_i s_i^2$;
that is, $\Phi> 0$.
Since $\Phi$-tensor describes the second fundamental form completely the
last property is equivalent to the fact that
$M\i \RR^q$ is \textit{stably} locally convex; that is, if $\Phi>0$, then $M$ is locally convex in $\RR^q\times \RR^k\supset\RR^q\times 0\z=\RR^q$ for some $k$.
(Given $k\in \NN$ there are examples of submanifolds $M\i \RR^q$
that are not convex as a submanifold in $\RR^q\times \RR^k$,
but convex as a submanifold in $\RR^q\times \RR^{k+1}$.)

\parbf{1.E. Positiveness of curvature tensor and symmetric $4$-tensors.}
The space $\S^2(\S^2(\T))$ splits into two subspaces,
the first is $\S^4(\T)\i \S^2(\S^2(\T))$ and the second is
$\A^4_+(\T)=\S^4(\T)^\bot\cap \S^2(\S^2(\T))$
which is canonically isomorphic to space of algebraic curvature tensors
$\A^4(\T)=\Lambda^4(\T)^\bot\cap \S^2(\Lambda^2(\T)$.

\medskip

\noindent
$\A^4(\T)$. 
Let $\mathcal{S}^*$ be the cone in $\A^4(\T)$ that consists of all tensors
$$\Rm(X,Y,Z,T)=\sum_i (s_i(X,Z)\cdot s_i(Y,T)-s_i(X,T)\cdot s_i(Y,Z)),$$ where
$s_i$  are positive elements of $\S^2(T)$. 
We say that such
curvature tensors have \emph{positive cosectional curvature};
this can be written as $\cosec (\Rm)>0$. 
For a Riemannian manifold $M$ we will write $\cosec (\Rm_p)>0$ or $\cosec_p>0$ if
the curvature tensor at $p\in M$ has positive cosectional curvature and
$\cosec M>0$ if the cosectional curvature of $M$ is positive at
all $p\in M$.

The curvature tensors with positive cosectional curvature are exactly
those that can be curvature tensors of submanifolds with positive
$\Phi$-curvature (equivalently, strictly convex submanifolds; see 1.D);
that is,
\begin{align*}
\cosec(&\Rm)>0
\\
&\Updownarrow
\\
\Rm(X,Y,Z,W)=\Phi(X,Z,Y,W)&-\Phi(X,W,Y,Z) \text{ for
some } \Phi>0.
\end{align*}
As you will see, any closed Riemannian
manifold $M$ with $\cosec M>0$ admits a smooth isometric embedding
$i\:M\to \RR^q$ with $\Phi(i)>0$.
In fact, we can assume that $q=\tfrac12\cdot(n+2)\cdot(n+5)$, see \cite[3.1.5(A) and 3.1.2(C)]{Grom-PDR}. 
The closure $\overline{\mathcal{S}}^*$ of  $\mathcal{S}^*$ can be
described as a minimal convex $\O(\T)$-invariant cone that contains the curvature tensor of product metric space $\SS^2\times \RR^{n-2}$.

The dual cone $\mathcal{S}$ (see 1.D) consists of curvature tensors with
positive sectional curvature.
For a point $p$ in a Riemannian manifold,
we will write $\sec_p>0$ or $\sec(\Rm_p)>0$ meaning
that $\Rm_p\in \mathcal{S}$.

The cones $\mathcal{S}^*$ and $\mathcal{S}$ are the smallest and largest
$\GL(\T)$-invariant cones in $\A^4(\T)$;
other $\GL(\T)$-invariant cones lie between $\mathcal{S}^*$ and $\mathcal{S}$.
The cone $\mathcal{Q}$ of
all curvature tensors with positive curvature operator is one of them;
it can be defined as
$$\mathcal{Q}
=
\set{R\in \A^4(\T)\i \S^2(\Lambda^2(\T))}{R=\sum_i\phi_i^2 \text{
for } \phi_i\in \Lambda^2(\T)}.$$
If the dimension is large, then $\mathcal{Q}\not=\mathcal{S}^*$.
Namely if dimension $=2$ or $3$,
then $\mathcal{Q}=\mathcal{S}^*$, and therefore $\mathcal{S}=\mathcal{Q}=\mathcal{S}^*$.
In dimension $4$, we have $\mathcal{Q}=\mathcal{S}^*$ and
$\mathcal{Q}^*=\mathcal{S}$.
The latter follows from \emph{Thorpe's characterization} of curvature tensors with positive
sectional curvature,  namely, if $M$ is positively curved
$4$-manifold, then there is a function $f$ on $M$ such that
$\Rm_x+f(x)\cdot \omega_x\z\in \S^2(\Lambda^2(\T_x))$ is a section of positive
quadratic forms on $\Lambda^2(\T)$; here $\omega$ denotes the
volume form, a section of $\Lambda^4(\T)\i \S^2(\Lambda^2(\T))$, see
\cite{Zol} for details.

In dimensions $5$ and higher, all the inclusions 
\[\mathcal{S}\subset \mathcal{Q}^*\subset\mathcal{Q}\subset\mathcal{S}^*\]
are strict.
Indeed: evidently, the inclusion $\mathcal{Q}^*\subset\mathcal{Q}$ is strict.
The inclusion $\mathcal{S}\subset \mathcal{Q}^*$ is strict if and only if so is
$\mathcal{Q}\subset\mathcal{S}^*$.
The latter is shown by example, see
\cite{Zol}.\footnote{In \cite{Grom-SGMC}, Gromov states the opposite. He writes \textit{ ``The closer of this cone (given by $Q\ge 0$)} [this is, $\overline{\mathcal{Q}}$ in our notations] \textit{can be defined as the minimal
closed convex $\O(n)$-invariant cone which contains the curvature
of the product metric on $\SS^2\times \RR^{n-2}$} [it is our $\overline{\mathcal{S}}^*$].''\dots}

\medskip

\noindent
$\S^4(\T)$.
Consider the cone of all forms $E\in \S^4(\T)$ such that
$E=\sum s_i^{\circ 2}$, where all $s_i$ are positive-definite and $s_i^{\circ 2}$ denotes the symmetric square of $s_i$.
It will be called the cone of \emph{positive forms} and denoted by $C_+$.
We will write $E>0$
if $E\in C_+\i \S^4(\T)$. Again, tensor $E$ is positive if it is a
symmetric part of a positive $\Phi$-tensor in
$\S^2(\S^2(\T))$.\footnote{By the way $C_+$ is also the
smallest $\GL(\T)$-invariant cone in $\S^4(\T)$. The biggest such cone
$C^*_+$ consists of all symmetric $4$-form $E$ such that
$$E(X,X,X,X)>0$$ for any nonzero $X\in \T$. 
Gromov in
\cite[3.1.4]{Grom-PDR} states that a symmetric form in
$E\in \S^{2k}(\T)$ is positive if and only if correspondent
quadratic form $E(\S^2(\T^k))\to \RR$ is positive-definite.
In our notations, it is equivalent to $C_+=Q_+$, and
this is equivalent to $C_+^*=Q_+^*$. 
The cone $C_+^*$ is the set of all positive-definite  forms in $\S^{2k}(\T)$;
equivalently, it is the set of positive homogeneous degree $2k$ polynomials on $\T$.
Analogously the cone $Q_+^*$ is the set of homogeneous degree $2k$ polynomials on $\T$ that can be expressed as a sum of squares.
Therefore this statement is equivalent to the following: \textit{each positive
polynomial is a sum of squares of polynomials}, and this was shown to be wrong in general.
Namely, David Hilbert \cite{Hil} had shown that this statement is true \textit{only} in the following three cases: (i)
$\dim \T\le 2$ and any $k$,
(ii) $k=1$ and any $\dim \T$,
(iii) $k=2$ and $\dim \T=3$.
This does not affect the rest of Gromov's book, except that one should always use the $C_+$-sense for positiveness.}

\section{Proofs}

I will not give here a proof of the first part of the local theorem for two reasons: first, it is real pain to write and read, and second, it is not mine.
The proof I have is a modification of Perelman's unpublished result.
Since this result was never published and (as far as I know) was never written,
I discuss it in Appendix \ref{appendix:B}.
(It should be more fun to look at the original proof than at my compilations.)

\parit{Proof of the second part of the local theorem.}
Let us prove first that
if $(M,g)$ is a Riemannian manifold with $\cosec M>0$,
then $(M,g)$ is isometric to a convex submanifold in $\RR^q$.
This is equivalent to the fact that there is an isometric embedding
$i\:M\to \RR^q$, such that
$\Phi(i)>0$ (see 1.D).

In general, a smooth isometric embedding of $(M,g)$ with $\cosec M>0$
may have an indefinite $\Phi$-tensor, but there is a way to make it positive.

Consider any smooth free isometric embedding  $i\:(M,g)\to\RR^q$.
By
Theorem \cite[3.1.5(A)]{Grom-PDR} for any tensor field $E\i \S^4(\T)$
such that $E_x>0$ (see 1.E) at all $x\in M$
one can find a $C^1$-close isometric embedding
$i'\:(M,g)\z\hookrightarrow \RR^q$, such that $E(i')=E(i)+E$.
In particular, one may choose $E=c\cdot g^{\circ 2}$ for any $c>0$.
Since $\cosec M>0$, for sufficiently large the $c$ we have $\Phi(i')>0$.

(Let us describe a more direct way to construct $i'$ with $E(i')=E(i)+c\cdot g^{\circ 2}$.
First construct an isometric embedding $j\:\RR^q\to \RR^{q'}$,
such that
$E(j)=c\cdot h^{\circ 2}$, where $h=\sum_{i=1}^q (dx_i)^2$
is the unit $2$-form on $\RR^q$ and then take $i'=j\circ i$.
One can construct $j$ on the following way:
first choose a collection of linear functions $l_i\:\RR^q\to \RR$
such that $h^{\circ 2} = \sum_i (dl_i)^4$,
second consider the diagonal
of the product of the following mappings: a linear mapping $L\:\RR^q\to \RR^q$
and twists $\tau_i\: \RR^q\to \RR^2$,
\[\tau_i(x)=
\bigl(a_i\cdot  \sin \bigl(b_i\cdot l_i(x)\bigr),a_i\cdot  \cos \bigl(b_i\cdot l_i(x)\bigr)\bigr)\]
with appropriately chosen $a_i$ and $b_i$.)

Now, since $M\i \RR^q$ is a convex submanifold, there is an open set $U\i M$ that is an intersection of locally convex hypersurfaces $F_i$ with obtuse angles between each pair of outward normals everywhere on $U$.
Approximate each $F_i$ as a convex polyhedral hypersurface $F_i^\varepsilon$ keeping the angles obtuse.
The intersection $U^\varepsilon$ of all $F_i^\varepsilon$ is a polyhedral submanifold,
and it has curvature $\ge 0$ (see 1.B).
Cutting subdomains from
$U^\varepsilon$ if necessary one gets the needed approximation. \qeds

\parit{Proof of the global theorem.} Let us first assume that
$\T M$ is stably trivial.
We will represent our submanifold $M$ as
an intersection of locally convex hypersurfaces with angles between any pair of outward normals $>\pi/2$
everywhere on $M$.
Once it is done, repeating the construction from the local theorem will finish the proof.

The existence of such representation is equivalent to the
existence of a smooth section of orthonormal bases $\{e_i\}$ in
$\N  M$ such that $s=\sum_i s_i e_i$ with all positive-definite 
$s_i\in \S^2(\T)$.

Consider a cover $U_k$, $k\in \{1,2,\dots,n\}$
of $M$ such that on each $U_k$ there is a smooth section of
orthonormal bases $\{e_{i,k}\}\i \N M$ with the above properties.
Since $\T M$ is stably trivial we can assume that $\N M$ is a trivial bundle.
Therefore we
can extend these bases to all $M$, and get $n$ bases
$\{e_{i,k}\}$ for all $\N M$.
Therefore at each point we have an isometric rotation $E_{k,k'}\in \O(q-m)$
that sends $\{e_{i,k}\}$ to $\{e_{i,k'}\}$.
Without loss of generality, we can assume that correspondent mapping
$E_{k,k'}\:M\to \O(q-m)$ is null-homotopic.

Let us choose a smooth partition of unity $u_k\:M\to [0,1]$; so,
$u_k|_{M\backslash U_k}\equiv 0$  for any $k$ and
$\sum_k u_k(x)\equiv 1$ for all $x\in M$.
At each point $x\in U_k\i M$ we have
$\Phi_x\equiv \sum_i s_{i,k}^2$.
Therefore for each $x\in M$ we have
$\Phi_x\equiv \sum_{i,k} u_k(x)\cdot  s_{i,k}^2$.
Consider $n\cdot\N M=\N_1 M\oplus \N_2 M\oplus \dots\oplus \N_n M$,
the sum of $n$ copies of the normal bundle,
take the basis $\{e_{i,k}\}$ for $\N_k$, and consider the subbundle $\N_\Delta(M)$,
that is spanned by
$(\sqrt{u_1} \cdot e_{i,1}, \sqrt{u_2} \cdot E_{12}e_{i,1},\dots, \sqrt{u_n} \cdot E_{1n}e_{i,1})$.
It is a trivial subbundle with a trivial orthogonal subbundle.
Therefore, there is a bundle isomorphism $i\:n\cdot \N M\z\to \RR^{(n-1)\cdot (q-m)}\times \N M$ that is an isometry on each fiber and sends
$\N_\Delta M$ to $\N M$.
Moreover, if $p_\Delta\:\N^n\to \N_\Delta$ is the orthogonal projection, then $i\circ p_\Delta (e_{i,k})\z=\sqrt{u_k} \cdot e_{i,k}$.

Therefore, we get a smooth section of orthonormal bases
$\{e_{i,k}\}\i \N'M\z=\RR^{(n-1)\cdot \dim(\N_x)}\times \N  M$;
that is, if we had $\N M$ as a normal bundle of $M\i \RR^{q}$, then $\N'M$ is a
normal bundle of $M\i \RR^{q}\times 0\i \RR^{q}\times
\RR^{(n-1)\cdot (q-m)}$.
Now for each pair of indexes $i,k$ we have a nonnegative quadratic form $s_{i,k}= \langle s,e_{i,k}\rangle$, $\Phi_x\equiv \sum_{i,k} u_k(x)\cdot s_{i,k}^2$, and at each point we have at least one quadratic form that is strictly positive.
It is not hard to rotate basis $e_{i,k}$ a little to get a new smooth section of bases in $\N'M$ with representation $s= \sum_{i,k}e_{i,k}\cdot s'_{i,k}$ where each $s'_{i,k}$ is strictly positive.

If $\T M$ is not stably trivial, then one still can find an embedding $M\hookrightarrow \RR^q$
that has positive $\Phi$-curvature at each point.
Take a small tubular neighborhood $U$ of $M$.
Let $\tilde M$ be a finite cover of $M$ such that $\T \tilde M$ is stably trivial.
 We can assume that $\N\tilde M$
is trivial, therefore $\N M$ is equivalent to a flat bundle.
From the above, we get the existence of a flat bundle $U'\to U$ such that the new induced normal bundle of $M$ with
respect to $U'$ is trivial.
The manifold $U'$ is flat, and repeating the same construction as above proves the theorem.
\qeds

\section{Remarks on curvature bounded above}

One can also ask a similar question for approximations by polyhedral spaces 
with upper curvature bound in the sense of Alexandrov; see \cite{AKP}.

\begin{thm}{Local theorem}
Let $P_n$ be a sequence of $m$-dimensional
polyhedral spaces
with curvature $\le \kappa$ that Lipschitz converges to a
Riemannian manifold $(M,g)$ of the same dimension, then $\cosec M\le \kappa$.
\end{thm}

Unfortunately, I cannot find as nice characterization for the global theorem
as in the case of curvature bounded below. Here is what I can do:

\begin{thm}{Global theorem}
Suppose $(M,g)$ is an $m$-dimensional Riemannian with $\cosec M\z\le \kappa$.
Assume that $M$ is diffeomorphic to a direct product of manifolds that care constant negative curvature, then $M$ can be realized as a Lipschitz limit of a sequence of $m$-dimensional polyhedral metrics with curvature at most $\kappa+\varepsilon$ for arbitrary $\varepsilon>0$.
\end{thm}

The proof of the local theorem is practically the same as for curvature bounded 
below. Perelman's lemma (which is the main technical tool in the proof) is 
also true for negative curvature.

The proof of the global theorem is also very similar, but one should consider 
embeddings into (noncomplete) flat pseudo-Riemannian manifold that is locally isometric to $\RR^{m,q}$ as a space-like submanifold of maximal dimension $m$ with a trivial normal bundle.
I do not have a complete answer to the question \textit{which manifolds admit such an embedding}.
If $M$ admits such an embedding, then the universal cover of $M$ is diffeomorphic to $\RR^m$. 
In particular, $M$ is $K(\pi,1)$ space, but so far I cannot say much useful 
about $\pi$.
On the other hand, if $M$ cares a metric of constant negative curvature, then it is a factor of pseudosphere in $\RR^{m,1}$ along some group,
and the factor of a little neighborhood of the pseudosphere by this group gives the needed ambient manifold.
Taking the product of such manifolds one has the needed embeddings for products of such manifolds.

\section{Problem section}\label{sec:Problem section}

The opposite question
\textit{which polyhedral metrics could be smoothed to a Riemannian manifold with positive curvature},
is wide open.
All examples I know so far satisfy the following conjecture:

\begin{thm}{Conjecture} Any polyhedral metric with curvature $\ge \kappa$
can be smoothed into a Riemannian orbifold with cosectional
curvature $\ge \kappa-\varepsilon$.
\end{thm}

The 2-dimensional case of the conjecture is a 
corollary of Alexandrov's embedding theorem \cite{Al}.
The 3-dimensional case is proved in \cite{LMPS}.
The conjecture would imply, in particular,
that any simply connected
manifold with positive cosectional curvature is diffeomorphic
to a sphere.
Indeed Corollary \ref{0.3}
implies that if $M$
is a Riemannian manifold with $\cosec M\ge 1$, then it can be approximated by polyhedral spaces $X_n$
with curvature $\ge 1$.
Therefore the spherical suspension $\SS(M)$ is approximated by $\Sigma (X_n)$ (cf. \cite{GW}).
From the conjecture, it would follow that $\Sigma(X_n)$ is smoothable into a Riemannian orbifold.
Therefore $M$ is a quotient of a sphere.

\medskip

Another question is whether the condition on stable triviality of the tangent bundle can be removed from Theorem \ref{0.2}.
 So far, I cannot even construct an approximation of
$(\CP^2,\text{can})$ by polyhedral metrics with curvature
$\ge -\varepsilon$.

One may also ask whether it is possible to construct an approximation
of $(\CP^2,\text{can})$ by polyhedral metrics with curvature
$\ge 0$. This is already a rigid question, in particular, from
Cheeger's results \cite{Ch} it is easy to see that any nonnegatively
curved polyhedral metric on $\CP^2$ carries a complex structure.
As it was pointed out by Mikhail Gromov,
$\CP^2$ carries polyhedral
metrics with curvature $\ge 0$. 
For example, if $P$ is a nonnegatively curved polyhedral that is homeomorphic to $\SS^2$, then its space of pairs $(P\times P)/\ZZ_2$ is homeomorphic to $\CP^2$ and naturally comes with a nonnegatively curved polyhedral metric.
This and many other examples are discussed by Dmitri Panov \cite{Pan}, but they do not solve our problem. 

\medskip

\textit{Can one generalize the Alexandrov embedding theorem?}
Namely, is it possible to characterize Riemannian manifolds that are isometric to a complete convex hypersurface in a complete convex hypersurface in \dots in $\RR^q$? 
Is it true that any simply
connected Riemannian manifold $M$ with $\cosec M>0$ is isometric to one of
those? 
If $M$ is compact it would imply that any such manifold is diffeomorphic to the standard sphere,
the latter follows from the result of Christoph Böhm and Burkhard Wilking \cite{BW}.

\medskip

Is it possible to characterize $m$-manifolds that admit embeddings into flat open $(m,q)$-pseudo-Riemannian manifold as a space-like surface?
(It is easy to see that if such an embedding exists, then 
the universal cover of $M$ is diffeomorphic to $\RR^m$, plus the first homotopy group must be linear, but I do not think it is a sufficient condition for the existence of such embedding.)

\appendix

\section[Example of positive curvature tensor with nonpositive cosectional curvature]{Example of positive curvature tensor with\\ nonpositive cosectional curvature}

Here I present calculations of Jost Eschenburg showing that
curvature tensor $R$ of $\SU(3)$ with bi-invariant metric has
nonnegative curvature operator but it is not true that $\cosec M\ge 0$.
This gives an example in dimensions $\ge 8$, from the work of
Zoltek \cite{Zol} it follows that such examples exist in dimensions $\ge 5$, but the calculations below are much simpler, and it might be useful if the reader wants quickly convince himself that such monsters do live.

\medskip

Consider the adjoint representation
$\ad\:\mathfrak{su}(3)\to \Lambda^2(\mathfrak{su}(3))$ for $\mathfrak{su}(3)$.
The curvature operator of $\SU(3)$
with bi-invariant metric has curvature operator
$R\:\Lambda^2(\mathfrak{su}(3))\to \Lambda^2(\mathfrak{su}(3))$ which coincides with
projection on $\Im(\ad)$.

Note that if $\Im(\ad)$ has no simple bi-vectors,
then the curvature operator of $\SU(3)$ with bi-invariant metric
does not have positive cosectional curvature.

Therefore we only have  to show that if $0\neq x\in \mathfrak{su}(3)$, then
$\ad_x\in \Lambda^2(\mathfrak{su}(3))$ is not a simple bi-vector; that is, $x\neq
v\wedge w$.

It is sufficient to prove it for $\ad_x$, where $x$ is tangent to a
maximal torus of diagonal elements in a matrix representation.
Therefore in the matrix representation it looks like
$x=\diag\{ai,bi,ci\}$ with $a+b+c=0$. Take the standard real basis
in $\mathfrak{su}(3)$ that comes from matrix form;
that is, take
 $A_1\z=\diag \{i,0,-i\},A_2=\diag \{0,i,-i\}$, take
 $F_1=e_2\wedge e_3,F_2= e_3\wedge e_1,F_3=e_1\wedge e_2$ be real and
$E_1=i\cdot e_2\circ e_3,E_2=i\cdot e_3\circ e_1,E_3=i\cdot e_1\circ e_2$
imaginary parts of basis, here $e_1,e_2,e_3$ is a basis of $\CC^3$
where $\SU(3)$ acts.
By the direct calculation we have
$\ad_x=(c-b)\cdot F_1\wedge E_1+(a-c)\cdot F_2\wedge E_2 + (b-a)\cdot  F_3\wedge E_3$,
now the fact that bi-vector $\phi\in \Lambda^2(\T)$
is simple is equivalent to $\phi\wedge\phi=0$, and
\begin{align*}
\ad_x\wedge \ad_x&=(c-b)\cdot(a-c)\cdot F_1\wedge E_1\wedge F_2\wedge E_2+
\\
&+(a-c)\cdot(b-a)\cdot F_2\wedge E_2\wedge F_3\wedge E_3+
\\
&+(b-a)\cdot(c-b)\cdot F_3\wedge E_3\wedge F_1\wedge E_1.
\end{align*}
Therefore if $\ad_x$ is simple, then at least two
of numbers $(c-b),(a-c),(b-a)$ are zeros
and since $a+b+c=0$ we have that  $a=b=c=0$; that is, $x=0$. \qeds

\section{Perelman's theorem and why I need it}\label{appendix:B}

In this appendix I will present the proof of one unpublished
result of G. Perelman;
it should be close to the original proof but some steps might differ.
It describes the main idea in the proof of the first part of the local theorem.

Let $M$ be an Alexandrov $m$-space and $U\i M$ be an open subset.
Let $F\:U\to \RR^m$ be a chart $F(p)=(x_1(p),x_2(p),\dots,x_m(p))$.
We say that $F$ is convex if each of the coordinate functions $x_i$
is convex.
The proof of the following claim easily follows from \cite[Proposition 3]{Per}

\begin{thm}{Claim}
Let $g$ be a convex function on $U$.
Suppose that for some
convex chart $F\:U\to \RR^m$ we have $\tfrac{\partial g}{\partial x_i}< 0$ for each $i$.
Then $g\circ F^{-1}$ is a convex function on $F(U)$. Moreover, for
any $p\in U$ and  $v\in \T_p$ we have
$$\nabla^2_vg\le\nabla^2_{dF(v)}\bigl(g\circ F^{-1}\bigr).$$
\end{thm}

In particular, if $S$ is a level surface of $g$, and $F$ is distance-nonexpanding, then 
$$I\!I_S(X,X)\le I\!I_{F(S)}(dF(X),dF(X))$$
for any $p\in S$.

\begin{thm}{Definition}
Let $X_n \buildrel {GH} \over \rightarrow X$ be a converging sequence of metric spaces and $f_n\:X_n\to X$ be a correspondent sequence of Hausdorff approximations.
We say that a sequence of measures $\mu_n$ on $X_n$ weakly converges to
measure $\mu$ on $X$ if for any continuous function $\alpha$ with compact
support on $X$ we have
\[\int_{X_n} \alpha\circ f_n\cdot \mu_n\to \int_{X} \alpha \cdot\mu.\]
\end{thm}

\begin{thm}{Theorem}\label{B.3.} Let $M_n$ be a sequence of Riemannian $m$-manifolds with curvature
$\ge \kappa$ that Lipschitz converges to a closed Riemannian manifold $M$.
Then scalar curvature on $M_n$ converges weakly to the scalar curvature on $M$;
that is, $\Sc_{g_n}\cdot \vol_{g_n}\rightharpoonup\Sc_{g}\cdot \vol_{g}$.
\end{thm}

If one has no lower bound for curvature, then the limit of scalar curvatures might be smaller than the scalar curvature of the limit \cite{Loh-SCH}.
It is unknown whether it could also be bigger.

The following lemma is a partial case of the theorem:

\begin{thm}{Lemma}
Let $F_n$ be a sequence of smooth convex hypersurfaces in $\RR^{m+1}$ that Hausdorff converges to a smooth convex hypersurface $F$.
Let $\Sc(F)$ and $\Sc(F_n)$ denote scalar curvatures of $F$ and $F_n$
and $h(F)$, $h(F_n)$ denote the $m$-Hausdorff measures of the
correspondent hypersurface.
Then $\Sc(F_n)\cdot h(F_n)$ converges weakly to $\Sc(F)\cdot h(F)$
\end{thm}

\parit{Proof.}
Let $\alpha$ be a continuous function with compact support in $\RR^{m+1}$.
Let us denote by $C_r(F)$ the set of points in $\RR^{m+1}$ that lie on outgoing normal rays to the hypersurface $F$ on the distance $< r$ to the hypersurface.
Let us define $\alpha_F\:C_\infty(F)\to \RR$ by $\alpha_F(x)=\alpha(y)$, where $y\in F$ is a closest point on the hypersurface.

It is well known and easy to see that
$\int_{C_r(F)}\alpha_F\cdot \vol$ is a polynomial of degree $m$ on $r$.
Moreover, its quadratic term is exactly $r^2\cdot\int_F\alpha\cdot \Sc(F)
\cdot h(F)$.

If $F_n\to F$, then $C_r(F_n)\to C_r(F)$
and $\alpha_{F_n}$ converge to $\alpha_{F}$.
Therefore, $\int_{C_r(F_n)}\alpha_{F_n}\cdot \vol\to\int_{C_r(F)}\alpha_F\cdot \vol$ and the coefficient with $r^2$ of correspondent polynomials also converge.
\qeds

\parit{Proof of \ref{B.3.}.} We first want to construct special
distance-like charts in a neighborhood of any point in $M$
together with some nice approximating charts on~$M_n$.

\begin{thm}{Lemma}\label{B.5} Given $p\in M$, $v\in \T_p^*(M)$, and $\varepsilon>0$ there is $\delta>0$,
a sequence $M_n\ni p_n\to p\in M$,
and sequence of convex functions $f_n\:B_\delta(p_n)\subset M_n\to \RR$ that converges to a convex function $f\:B_\delta(p)\subset M\to\RR$ such that $d_pf=v$, $|f''|<\varepsilon$ everywhere on $B_\delta(p)$.
\end{thm}

\parit{Proof of \ref{B.5}.}
Consider an orthonormal basis $\{e_i\}$
in $\T_p M$ such that $\sum_i e_i=c\cdot v$. Take $r >0$ an let
$a_i=\exp_p (r\cdot e_i)$.
Let $f=\sum_i \phi\circ\text{dist}_{a_i}$,
where $\phi(x)=\alpha\cdot  \log x-\beta\cdot  x^2$ if $\dim M=2$ and
$\phi(x)=\alpha\cdot  \frac{1}{x^{n-2}}-\beta\cdot  x^2$ if $\dim M>2$.
The same arguments as in \cite[4.3]{PP} shows that for appropriately chosen constants $\alpha$ and $\beta$, the function $f$ satisfies the conditions
of the lemma in a small ball $B_\delta(p)$ .

To construct an approximation of $f$, construct a
sequence $a_{i,n}\to a_i$ for each $i$, and take $f_n=\sum_i
\phi\circ\text{dist}_{a_{i,n}}$. Again the same reasoning as in
\cite[4.3]{PP} proves that there is $\varepsilon>0$ such that for large $n$
the function $f_n$ is convex in a $\delta$-neighborhood of $p_n$.\qeds

Choose any orthonormal basis $v_1,\dots,v_m\in \T^*_p M$.
For each $v_i$, use the lemma to construct function $f_i$ on $B_\delta(p)$ together with its approximations
$f_{i,n}$ on $B_\delta(p_n)$.
In addition to the above properties, these functions will be almost orthogonal for a small $\delta>0$;
that is, one can assume that angle between level surfaces lies in the range $\tfrac\pi2\pm\varepsilon$.

Now we start induction by dimension, we can take $\dim =2$ as a base, in which case convergence follows from the Gauss--Bonnet formula.
Assume we already proved it for all dimensions $<m$.

To save space and time, let us agree that extra index
$n$ will always denote correspondent babe for $M_n$.

Let $p\in M$ and $S_1,S_2,\dots,S_m$ be one-parameter families of
 coordinate surfaces $f_i=c$.
Let us denote by $\Sc_i$ is the ``scalar'' curvature of directions tangent  to $S_i$; in other words, $\Sc_i=\Sc-\Ricc(u_i)$ where $u_i$ is the unit vector field normal to $S_i$.
Note that from the lower curvature bound we have $|\text{Rm}|<c_1+c_2\cdot \Sc$
and therefore $(1\pm\alpha)\cdot(m-2)\cdot \Sc\gtrless\Sc_1+\Sc_2+\dots+\Sc_m$
where $\alpha$ depends on angles between these coordinate
surfaces and $\alpha\to 0$
as all these angles converge to $\pi/2$, in particular, as $\varepsilon\to 0$.

Let $\Sc(S_i)$ be the scalar curvature of the intrinsic metric of the correspondent coordinate surface.
Since the Jacobian of our charts converges to the Jacobian of the limit chart, from the induction hypothesis we have $\Sc({S_{i,n}})\cdot \vol_{g_n}$
converges weakly to $\Sc(S_i)\cdot \vol_g$.

Denote by $\G(S_i)$ the Gauss curvature of $S_i$;
that is, $\G(S_i)\z=\sum_{j'\not=j}k_{j'}\cdot k_j$, where $k_j$ are the principal curvatures of $S_i$.
By the Gauss formula, we have $\Sc_i+\G(S_i)= \Sc(S_i)$.
Since each $S_i$ is convex
$$\Sc_i\le \Sc(S_i)\ge \G(S_i).$$
Therefore, after passing to a subsequence,
$\Sc_{i,n}$ should converge weakly to some
$\overline{\Sc}_i\le \Sc(S_i)\le \Sc_i+n\cdot (n-1)\cdot \varepsilon^2$.

Let us prove the following lower bound: $\overline {\Sc}_i\ge  \Sc_i-C\cdot\varepsilon^2$ for some $C=C(m)$.
By the inequality above, after passing to a subsequence $\G(S_{i,n})$
 converges to some measure $\overline \G(S_i)$ and obviously,
\[\overline \Sc_i=\Sc(S_i)-\overline \G(S_i)\ge \Sc_i-\overline \G(S_i).\]
Therefore it is sufficient to show that $\overline \G(S_i)\le C\cdot\varepsilon^2$ for some fixed $C$.

To prove this last estimate, let us construct a new chart similar to
one before, $H=(h_1,h_2,\dots,h_m)$ with the approximations
$H_n=(h_{1,n},h_{2,n},\dots,h_{m,n})$ such that $\tfrac{\partial f}{\partial
h_i}<\tfrac{-1}{10\cdot m}$.
From the claim above we have that $\G(S_n)\le c
\cdot\G(H(S_n))$ as well as $\G(S)\le c\cdot \G(H(S))\le C\cdot\varepsilon^2$

Now $H(S_n)$ converges to $H(S)$ as convex hypersurfaces in $\RR^m$.
Applying the lemma, we get $\G(H(S_n))$ converges weakly to $\G(H(S))$.

The theorem follows since for any $\varepsilon>0$ there is a finite covering of $M$ by charts as in \ref{B.5}.
\qeds

Along the same lines one can prove stronger statements; see \cite{LP} for details.

\begin{thm}{Smooth Proposition} Let $(M_n,g_n)$ be a sequence of
Riemannian $m$-manifolds with curvature $\ge \kappa$ that
GH-converges to a Riemannian manifold $(M,g)$ of the same
dimension $m$.
Then there is a sequence of reparameterizations (diffeomorphisms) $f_n\:M\to M_n$, such that the curvature tensor of $df_n^*(g_n)$ weakly converges to the curvature tensor of $g$ on $M$.
\end{thm}

\begin{thm}{Corollary} Let $R$ be an $\SO(\T)$ invariant convex set in $\A^4(\T)$. Assume
that there is a lower bound $\kappa>-\infty$ for sectional curvature in $R$.
Let $M_n$ be a sequence of Riemannian manifolds with curvature tensor from $R$ at each point.
Suppose $M_n$ converges to a Riemannian manifold $M$ of the same dimension.
Then the curvature tensor at any point of $M$ is from $R$.
\end{thm}

For example, a smooth limit of manifolds with positive curvature operator must have positive curvature operator.

Note that this corollary cannot hold for general
$\SO(\T)$ invariant convex set in $\A^4(\T)$, for example as it shown in \cite{Loh-GLC}, \cite{Loh-SCH} it is not true for sets $R=\{r\z\in
\A^4\:\Ricc(r)\le c\}$ curvature and for $R'=\{r\in \A^4\mid c\le
\Sc(r)\le c+\varepsilon\}$.
But, I believe that the condition on lower sectional curvature $R$ can be relaxed.

Finally, one can give a singular version of this result that we need in our paper:

First, let us describe the singular curvature tensor of a polyhedral space. 
Assume
we have a $(1\pm\varepsilon)$-bi-Lipschitz parametrization of polyhedral $P$
by smooth
Riemannian manifold $f\:M\to P$, such that $f^{-1}$ is smooth on each simplex.
One can think about $P$ as $(M,d)$ where $d$ is a singular metric.
Now let us define the curvature tensor of $d$ as follows: its support is the image of the $(n-2)$-skeleton of $P$ and on each $(n-2)$ simplex it is defined as $h_{n-2}\cdot (2\pi-\omega)\cdot\alpha^2$ where $h_{n-2}$ is the Hausdorff
measure of this image, $\omega$ is the total angle around this simplex $\Delta$ and $\alpha=dx\wedge dy$ is a bi-vector field with the following properties: $|\alpha|=1$ everywhere on the image of simplex and $\alpha|_{f^{-1}(\Delta)}=0$.

\begin{thm}{Singular Proposition} Let $P_n$ be a polyhedral $m$-spaces
with curvature $\ge \kappa$ that GH-converges to a Riemannian
manifold $(M,g)$ of the same dimension $=m$. Then there is a
sequence of smooth parameterizations $f_n\:M\to P_n$, such that the
described singular curvature tensor  weakly converges to the
curvature tensor of $g$ on $M$.
\end{thm}

As well as in the corollary above, since the cosectional curvature of $(M,d_n)$ is positive, we get that the curvature tensor on the limit $(M,g)$ has positive cosectional curvature.
The latter implies the first part of the local theorem (\ref{0.1}).

\sloppy
\def\emph{\textit}
\printbibliography[heading=bibintoc]
\fussy

\end{document}